\documentclass{article}
\usepackage{graphicx} \graphicspath{ {seipati/} }
\usepackage{amsmath}
\usepackage{color}
\usepackage{fancyhdr}
\usepackage{verbatim}
\usepackage{amsfonts,longtable}
\usepackage{epsfig}
\usepackage{amsfonts}
\usepackage{color}

\def\R{\hbox{{\rm I}\kern-0.2em{\rm R}\kern0.2em}}

\def\d{{\rm d}}

\def\bn{\begin{equation}}
\def\en{\end{equation}}
\def\bny{\begin{eqnarray*}}
\def\eny{\end{eqnarray*}}
\def\be{\begin{eqnarray}}
\def\ee{\end{eqnarray}}
\def\bc{\begin{center}}
\def\ec{\end{center}}
\def\X{{\cal X}} 
\def\p{\partial} 
\def\({\left(}
\def\){\right  )}
\def\[{\left[}
\def\]{\right]}
\def\bc{\begin{center}}
\def\ec{\end{center}}

\newtheorem{dfn}{Definition}[section]
\newtheorem{thm}{Theorem}[section]
\newtheorem{rem}{Remark}[section]
\newtheorem{pro}{Proposition}[section]

\newtheorem{cor}{Corollary}[section]
\newtheorem{lem}{Lemma}[section]
\newtheorem{exm}{Example}[section]

\def\bn{\begin{equation}}
\def\en{\end{equation}}
\def\bny{\begin{eqnarray}}
\def\eny{\end{eqnarray}}
\def\be{\begin{eqnarray*}}
\def\ee{\end{eqnarray*}}
\def\bdn{\begin{dfn}}
\def\edn{\end{dfn}}
\def\btm{\begin{thm}}
\def\etm{\end{thm}}
\def\bpf{\begin{proof}}
\def\epf{\end{proof}}
\def\bpn{\begin{pro}}
\def\epn{\end{pro}}
\def\brk{\begin{rem}}
\def\erk{\end{rem}}
\def\bcy{\begin{cor}}
\def\ecy{\end{cor}}
\def\blm{\begin{lem}}\def\elm{\end{lem}}
\def\bex{\begin{exm}}
\def\eex{\end{exm}}

\def\X{{\cal X}} 
\def\Y{{\cal Y}}  
 \def\R{{\cal R}}
 \def\L{{\bf L}} \def\Z{{\cal Z}}

\begin{document}

\title{\textbf{On the relationship between the invariance and conservation laws of differential equations}}
\author{A H Kara$^{*}$
 \\
\\
 School of Mathematics, University of the Witwatersrand,\\ Private Bag 3, Wits 2050,
Johannesburg, South Africa
}

\date{}
\maketitle
\begin{quote}

{\textbf{Abstract:} {\small In this paper, we highlight the complimentary nature of the results of Anco \& Bluman  and Ibragimov  in the construction of conservation laws; that whilst the former establishes the role of multipliers, the latter presents a formal procedure to determine the flows. Secondly, we show that there is an underlying  relationship between the symmetries and conservation laws in a general setting - extending the results of Kara \& Mahomed. The results take apparently differently forms for point symmetry generators and higher-order symmetries. Similarities exist, to some extent, with a previously established result relating symmetries and  multipliers of a differential equation. A number of examples are presented.
}}

 {\textbf{Keywords:} } {\small  conservation laws; symmetries  }

\end{quote}

\vfill

$^{*}$ Corresponding Author, E-mail: Abdul.Kara@wits.ac.za

\newpage

{\section{Introduction}}

The role and methods associated with conservation laws are now well established and there has been some momentous works in these areas in recent times
building on the contributions made by Noether which generally dealt with variational problems those that admit variational symmetries. It is not surprising then that much of the recent works focussed on generalizations as far as constructions of conservation laws go; possibly non variational and preferably independent of a knowledge of symmetries.
A vast amount and extensively cited works  are due to Anco \& Bluman in \cite{bgi,bgi2}, inter alia, Anderson \cite{iand,ian}, Kara \& Mahomed \cite{ukm} and a useful in depth treatise is presented in the work of Olver \cite{o} which goes a long way in discussing the concept of `recursion operators'. The first of these deals extensively with the notion of `multipliers'; that if a differential equation times a factor (differential function) is closed, then the Euler operator annihilates this product so that finding conserved flows amounts to finding the factors. It turns out that the multipliers are solutions of the adjoint equation. Of course, one still needs to determine the corresponding conserved flows using, amongst others, homotopy formulae \cite{h}. A large amount of software to construct the various components of conserved vectors are available, see \cite{h1,chev}.

Since conservation laws seem to be tied in with invariance properties, the intention to avoid the symmetry route can prove to be difficult. This is partly
due to the amount of work required to construct conserved flows; it can be cumbersome and tedious when dealing with the large systems of differential equations that arise in physics, cosmology and engineering. For e.g., constructing conservation directly from the definition may be straightforward for simple scalar ordinary differential equations but the more complex the differential equation, as they are in fluids, cosmology and the various systems of Schr\"odinger equations that is abundant in the literature (to name a few), the greater the task. The popularity of Noether's theorem lies in the existence of a formula. Trying to mimic this formula even in the non variational case has been tempting and partly successful, see \cite{partial}. In particular, the recent work of Ibragimov \cite{i} develops a procedure to construct conserved vectors using the Noether operator, a symmetry of the differential equation solutions of the adjoint equation.

An in depth study into the results due to Anco \& Bluman in \cite{bgi,bgi2} and Ibragimov \cite{i} suggests that similarities are abundant - see \cite{anco17}. However, it also shows that since the methods employed are largely different, there are some intrinsics differences and what is presented here is an attempt to show that these differences, in fact, allows these works to complement each other. For example, the underlying aspect in the multiplier approach is primarily to construct multipliers that leads to the  differential equation being conserved. These multipliers can be chosen with a specific order (in derivatives) in mind and then one may choose from a number of methods to construct the conserved vectors. In \cite{i}, the particular method appeals to the Noether operator after having knowledge of a symmetry and a solution of the adjoint equation. It will be shown, that the total divergence of the conserved flow has a form dependent on whether the symmetry used is a point symmetry or an evolutionary/canonical symmetry; the general result in the latter case would include generalised symmetries.

{\section{Notations and preliminaries}}

What follows is a summary of the definitions, concepts and notations that will be utilised in the sequel.

Consider an $k$th-order system of partial differential equations (pdes)
of $n$ independent variables $\underbar{x}=(x_1,x_2,\ldots,x_n)$ and $p$
dependent variables $w=(w_1,w_2,\ldots,w_p) $ viz.,
\begin{equation}
	E(\underbar{x},w,w_{(1)},\ldots,w_{(r)})= 0,~~~ ~~~ ~~~u =1,\ldots,\tilde{p},
\end{equation}
where a locally analytic function $f(\underbar{s}, w, w_1,\dots,w_k)$ of a finite number of dependent variables $w, w_1,\dots,w_b$ denote the collections of all first , second ,$\dots$, $b$th-order partial derivatives and $s$ is a multivariable, that is
\begin{equation}
w_i^\alpha = D_i(w^\alpha), ~~~ ~~~ ~~~ w_{ij}^\alpha = D_jD_i(w^\alpha),\dots
\end{equation}
respectively, with the total differentiation operator with respect to $x^i$ given by,
\begin{equation}
D_i = \frac{\partial}{\partial x^i} + w_i^\alpha \frac{\partial}{\partial w^\alpha} + w_{ij}^\alpha \frac{\partial}{\partial w_j^\alpha}+\dots ~~~ ~~~ ~~~ i = 1,\dots ,p.
\end{equation}

In order to determine conserved densities and fluxes, we resort to
the invariance and multiplier approach based on the well known
result that the Euler-Lagrange operator annihilates a total
divergence. Firstly, if $(T^{x1}, T^{x2},\ldots)$ is a conserved
vector corresponding to a conservation law, then
\begin{equation} D_{x1} T^{x1}+D_{x2} T^{x2}+\ldots=0\label{n1}\end{equation}
along the solutions of the differential equation $E(\underbar{x},w,w_{(1)},\ldots,w_{(k)})= 0$.\\
\\
Moreover, if there exists a nontrivial differential function $Q$,
called a `multiplier', such that
\begin{equation}\label{3.5a}
	 Q (\underbar{x},w, w_{(1)} \dots) E(\underbar{x},w,w_{(1)},\ldots,w_{(r)}) =D_{x1} T^{x1}+D_{x2} T^{x2}+\ldots ,
\end{equation}
for some (conserved) vector $(T^{x1}, T^{x2},\ldots)$, then
\begin{equation}\label{3.5b}
	{{\delta}\over {\delta u}}[ Q (\underbar{x}, w, w_{(1)} \dots) E(\underbar{s},w,w_{(1)},\ldots,w_{(r)})] = 0,
\end{equation}
where ${{\delta}\over {\delta w}}$ is the Euler operator. Hence, one may determine the multipliers, using (\ref{3.5b}) and then construct the corresponding conserved vectors; several approaches for this exists of
which the better known one is the `homotopy' approach.

If a pde is variational, then the conservation laws may be constructed from Noether's Theorem.  It can be shown that Lie point symmetries that leave the system of differential equations invariant contain the algebra of Noether/divergent/variational symmetries \cite{o,i}.

Conservation laws may be expressed as conserved forms \cite{ian}. For example, if $\underbar{x}=(t,s)$, the conserved form would be
$$\omega=T^t\d s-T^s\d t$$ (where $(T^t,T^s)$ is the conserved vector such that $D_tT^t+D_sT^s=0$ on the solutions of the pde $E(s,t,w,w_{(1)},\ldots,w_{(r)})= 0$ ). Here, $T^t\d s$ leads to the `conserved density' if $t$ and $s$ are time and space, respectively.

{\section{Conservation laws}}

\textbf{3.1} In the first case, we consider the relationship between the conserved flows and the respective point symmetry generators of the differential equation.

Example 1.

We, firstly, utilise the heat equation $u_t-u_{xx}=0$ as an illustrative example. The final result is presented in a Proposition. The multipliers $Q_1=-x$ and $Q_2=-e^t\sin x$ are discussed in \cite{i} to construct conserved vectors, there referred to as solutions of the adjoint equation $v_t+v_{xx}=0$. Thus, ${{\delta}\over {\delta u}}[Q_1(u_t-u_{xx} )]=0$. In general, then, ${{\delta}\over {\delta u}}[v(u_t-u_{xx} )]=0$ so that, by \cite{i}, the conserved flow via the point symmetry $X=2t\p_x-xu\p_u$
is
\bn\begin{array}{ll}
T^t&=-v(xu+2tu_x)   ,\\
T^x&= v(2tu_t+u+xu_x)-v_x(xu+2tu_x)
\end{array}\label{a}\en
The total divergence is
\bn\begin{array}{ll}
D_tT^t+D_xT^x&=-v_t(xu+2tu_x)-v(xu_t+2u_x+2tu_{xt})  +v_x(2tu_t+u+xu_x)\\
&+v(2tu_{xt}+2u_x+xu_{xx})-v_x(u+xu+x+2tu_{xx})-(xu+2tu_x)v_{xx}\\
&=-(xu+2tu_x)(v_t+v_{xx})-(2tv_x-xv)(u_t-u_{xx})
\end{array}\label{a1}\en
That is the total divergence takes the form
\bn\begin{array}{ll}
D_tT^t+D_xT^x&=-(xu+2tu_x)(v_t+v_{xx})+(-xv+2tv_x)  (u_t-u_{xx})\\
&=W(v_t+v_{xx})+(\X v)(u_t-u_{xx})\\
&=(\X v)(u_t-u_{xx})
\end{array}\label{a3}\en
where $\X=2tD_x-x$. If $v=1$ and $W=-(xu+2tu_x)$ is the characteristic of $X$. If, in $Q=-xv+2tv_x$  ,  \\
(i) $v=1$, then $D_tT^t+D_xT^x=-x(u_t-u_{xx})$ which leads to the multiplier  $Q_1=-x$ \\
(ii) $v=x$, then $D_tT^t+D_xT^x=(-x^2+2t)(u_t-u_{xx})$ leading to the multiplier $Q=-x^2+2t$.

Example 2.

Consider the one-dimensional wave equation $u_{tt}-u_{xx}=0$ and the Lorentz rotation symmetry $Y=t\p_x+x\p_t$ with characteristic $W=-tu_x-xu_t$ and adjoint equation $-v_{tt}+v_{xx}=0$. The detailed calculation using the results in \cite{i} leads to
\bn\begin{array}{ll}
T^t&=xv(u_{tt}-u_{xx})+v_t(xu_t+tu_x)+v(-xu_{tt}-u_x-tu_{xt}  ,\\
T^x&=tv(u_{tt}-u_{xx})-v_x(xu_t+tu_x)-v(-tu_{xx}-u_t-xu_{xt}
\end{array}\label{b}\en
so that
\bn\begin{array}{ll}
D_tT^t+D_xT^x&=-(x u_t+tu_x)(-v_{tt}+v_{xx})+(xv_t+t v_x)  (u_{tt}-u_{xx})\\
&=W(-v_{tt}+v_{xx})+(\Y v) (u_{tt}-u_{xx})\\
&=(\Y v) (u_{tt}-u_{xx}),
\end{array}\label{b3}\en
where $\Y =tD_x+xD_t$.

The following Proposition that defines the relationship between point symmetries, multipliers and conservation laws constructed via the Noether operator in \cite{i}, can be easily proved.

\textbf{Proposition 1}. If $Z=\xi(x,t,u)\p_x+\tau(x,t,u)\p_t+\phi(x,t,u)\p_u$ (characteristic $W=\phi-\xi u_x-\tau u_t $) is a Lie point symmetry generator of a second-order partial differential equation (pde) $E(x,t,u,u_{(x)},u_t,u_{tt},\ldots)= 0$ (whose adjoint equation is $F(x,t,v,v_{(x)},v_t,\ldots)= 0$), $\L=v E$ and
\bn\begin{array}{ll}
T^t&=\tau\L+W({{\p \L}\over{\p u_t}}-D_t{{\p \L}\over{\p u_{tt}}}  -D_x {{\p \L}\over{\p u_{tx}}} ) +D_tW {{\p \L}\over{\p u_{tt}}}+D_x{{\p \L}\over{\p u_{tx}}} ,\\         T^x&=\xi\L+W({{\p \L}\over{\p u_x}}-D_t{{\p \L}\over{\p u_{xt}}}  -D_x {{\p \L}\over{\p u_{xx}}} ) +D_tW {{\p \L}\over{\p u_{xt}}}+D_x{{\p \L}\over{\p u_{xx}}} ,\\
\end{array}\label{c}\en
then the divergence
\bn\begin{array}{ll}
&D_tT^t+D_xT^x\\
&=W F(x,t,v,v_{(x)},v_t,\ldots)+(\xi(x,t,v)v_x+\tau(x,t,v)v_t+\phi(x,t,v)    ) E(x,t,u,u_{(x)},u_t,u_{tt},\ldots)\\
&=W F(x,t,v,v_{(x)},v_t,\ldots)+(\Z v-\lambda v)  E(x,t,u,u_{(x)},u_t,u_{tt},\ldots),
\end{array}\label{c1}\en
where $\Z=\xi(x,t,v)D_x+\tau(x,t,v)D_t+\phi(x,t,v)$ and $\lambda$ is determined by the conformal factor. That is, if $ZE=\mu_1 E$ and $D_t\tau+D_x\xi=\mu_2$, then $\lambda=\mu_1+\mu_2$; $\lambda$ need not be a constant. On particular solutions $v=v(x,t)$ of the adjoint equation, we have a conserved flow $(T^t,T^x)$ with multiplier $Q=\xi(x,t,v)v_x+\tau(x,t,v)v_t+\phi(x,t,v(x,t))$.

After some cumbersome calculations, Proposition 1 is easily generalised to the multidimensional pde $E(x,t,u,u_{(x)},u_t,\ldots,u_{(r)})= 0$.

Example 3a.

For the the third order KdV equation $u_t-uu_x-u_{xxx}=0$, the adjoint equation is $v_t-v_{xxx}-uv_x=0$. The calculations in \cite{i}, using the point symmetry $X=-3t\p_t-x\p_x+2u\p_u$ and an extended version of (\ref{c}), lead to the conserved vector components
\bn\begin{array}{ll}
T^t&=v(3tu_{xxx}+3tuu_x+xu_x+2u  ) ,\\
  T^x&=-v(2u^2+xu_t+3tuu_t+4u_{xx}+3tu_{txx} )+v_x( 3u_x+3tu_{tx} +xu_{xx})-v_{xx}(2u+3tu_t+xu_x  ) ,\\
\end{array}\label{d}\en
so that, after detailed simplification we get
\bn\begin{array}{ll}
D_tT^t+D_xT^x&=(2u+3tu_u+xu_x  ) (v_t-v_{xxx}-uv_x  )+(-xv_x-3tv_t+2v  -v    ) (u_t-uu_x-u_{xxx} )\\
&=W (v_t-v_{xxx}-uv_x  )+[(\X  -1    )v] (u_t-uu_x-u_{xxx} ),
\end{array}\label{d1}\en
where $W=2u+3tu_u+xu_x$ and $\X=-3t\p_t-x\p_x+2$. In Proposition 1, $\lambda=5-4=1$.

Example 3b.

Consider the simplest Schr\"odinger equation with cubic nonlinearity $iu_t-u_{xx}+u|u|^2=0$. If we put $u=p+iq$ then define
$\textbf{L}=v[-q_t-p_{xx}+p(p^2+q^2)]+w[p_t-q_{xx}+q(p^2+q^2)]$ where $(v,w)$ is the solution of the system
$-w_t-v_{xx}+v(p^2+q^2)+2p(vp+wq)=0$, $v_t-w_{xx}+w(p^2+q^2)+2q(vp+wq) =0  $ - the adjoint of
$-q_t-p_{xx}+p(p^2+q^2)=0$, $p_t-q_{xx}+q(p^2+q^2)=0$. The components of the conserved vector, using $X=\p_t$ are then
\bn\begin{array}{ll}
T^t&={\bf L} -p_tw+vq_t,\\
  T^x&=-p_tv_x+p_{xt}v-q_tw_x+wq_{xt}\\
\end{array}\label{d3}\en
so that, after some manipulation,
\bn\begin{array}{ll}
&D_tT^t+D_xT^x\\
&=p_t[-w_t-v_{xx}+v(p^2+q^2)+2p(vp+wq)  ]  +q_t[v_t-w_{xx}+w(p^2+q^2)+2q(vp+wq)  ]\\
&+v_t[-q_t-p_{xx}+p(p^2+q^2)]+w_t[p_t-q_{xx}+q(p^2+q^2) ]\\
&=-W^1[-w_t-v_{xx}+v(p^2+q^2)+2p(vp+wq)  ]-W^2[v_t-w_{xx}+w(p^2+q^2)+2q(vp+wq)  ]\\
&-(\X v)[-q_t-p_{xx}+p(p^2+q^2)]-(\X w)[p_t-q_{xx}+q(p^2+q^2) ],
\end{array}\label{d4}\en
where $W^1=-p_t$, $W^2=-q_t$ and $\X=-D_t$. When $v=p$ and $w=q$, we get the well known energy conservation via (\ref{d3}) using the multiplier $(p_t,q_t)$.

\textbf{3.2} We now consider the connection between generalised symmetries, higher-order symmetries and evolutionary/canonical symmetries and associated conservation laws. Again,we suppose $\textbf{L}=v(x,t)E$ \cite{i}.

Example 4. In this example, we revisit the heat equation $u_t-u_{xx}=0$ with its evolutionary symmetry $X_1=(tu_x+\frac12 x u)\p_u$ (from the point symmetry $-t\p_x+\frac12 x u\p_u )$, higher symmetries $X_2=u_{xx}\p_u$ and $X_3=(2tu_{xxx}+xu_{xx})\p_u$ used to construct conserved flows $(T^t,T^x)$.

(i) With $X_1$, we obtain the components of the conserved vector to be
\bn\begin{array}{ll}
T^t&=v(tu_x+\frac12 x u ) ,\\
  T^x&=-v( tu_{xx}+\frac12  u+\frac12 x u_x ) ,\\
\end{array}\label{f}\en
so that
\bn\begin{array}{ll}
D_tT^t+D_xT^x&=v_t(tu_x+\frac12 x u )+v( u_x+tu_{xt}+\frac12 x u_t  )\\
&+-v_x( tu_{xx}+\frac12  u+\frac12 x u_x )-v(tu_{xxx}+u_x+\frac12 x u_{xx}   )\\
&=(v_t+v_{xx})(tu_x+\frac12 x u)+v(tu_{xt}-tu_{xxx}+\frac12x u_t\frac12xu_{xx}      )\\
&=W (v_t+v_{xx}  )+v\R_1  (u_t-u_{xx} ),
\end{array}\label{f1}\en
where $W=tu_x+\frac12 x u$ and $\R_1= tD_x+\frac12x$ is the recursion operator associated with $X_1$.

(ii) Using $X_2$, we get
\bn\begin{array}{ll}
T^t&=vu_{xx} ,\qquad  T^x=v_x u_{xx}-vu_{xxx} ,\\
\end{array}\label{f2}\en
so that
\bn\begin{array}{ll}
D_tT^t+D_xT^x&=u_{xx}(v_t+v_{xx}  )  +vD_xD_x(u_t-u_{xx})=W(v_t+v_{xx}  )  +v \R_2(u_t-u_{xx}).
\end{array}\label{f3}\en

(iii) With $X_3$, we get
\bn\begin{array}{ll}
T^t&=v(xu_{xx}+2tu_{xxx} ) ,\\
  T^x&=v_x(xu_{xx}+2tu_{xxx}  )-v(u_{xx}+xu_{xxx}+2tu_{xxxx}  ) ,\\
\end{array}\label{f4}\en
so that, after some simplifications the total divergence is
\bn\begin{array}{ll}
D_tT^t+D_xT^x&=(v_t+v_{xx})(xu_{xx}+2tu_{xxx} )+v( 2t[u_{txxx}-u{xxxxx}   ]+x[u_{txx}-u{xxxx}   ]        )\\
&=W(v_t+v_{xx}  )+v\R_3  (u_t-u_{xx} ),
\end{array}\label{f1}\en
where $W=2tu_{xxx}+xu_{xx}$ and $\R_3=xD_xD_x+2tD_xD_xD_x$ is the respective recursion operator.

Example 5. It is well known that, for the wave equation $u_{tt}-u_{xx}=0$ and any variational equation,  `multipliers' or, equivalently, solutions of the adjoint equation are symmetries of the equation so that, in the simple case of the evolutionary vector field $Y=(tu_x+xu_t)\p_u$ is a generalised symmetry and $Q=tu_x+xu_t$ is multiplier. Applying Noether's theorem is clearly the efficient route to constructing a conservation law. Alternatively, if we assume $L=v(x,t)(u_{tt}-u_{xx}  )$ using the procedure in \cite{i}, we get
\bn\begin{array}{ll}
T^t&=v_x(tu_x+xu_t  )  - v(tu_{xx}+u_t+xu_{xt} ) ,\\
  T^x&=-v_t(tu_x+xu_t  )+v(u_{x}+xu_{tt}+tu_{xt}  ) ,\\
\end{array}\label{g}\en
so that
\bn\begin{array}{ll}
D_tT^t+D_xT^x&=(-v_{tt}+v_{xx})(tu_x+xu_t )+v( t[u_{ttx}-u{xxx}   ]+x[u_{ttt}-u{txx}   ]        )\\
&=W(-v_{tt}+v_{xx} )+v\R  (u_{tt}-u_{xx} ),
\end{array}\label{g1}\en
where $W=tu_x+xu_t$ and $\R=tD_x+xD_t$.

\textbf{Proposition 2}. In Proposition 1, if $Z$ is a generalised symmetry or evolutionary/canonical vector field such that $ZE=(\R+\lambda)E$, where $\R$ is the recursion operator associated with $Z$, then
\bn\begin{array}{ll}
D_tT^t+D_xT^x&= W F(x,t,v,v_{(x)},v_t,\ldots)+v(\R+\lambda)E(x,t,u,u_{(x)},u_t,u_{tt},\ldots).
\end{array}\label{e1}\en

Again, the proposition can be generalised to the multi-dimensional case.

Example 6. We revisit the KdV equation with its evolutionary vector field $X=(xu_x+3tu_t+2u)\p_u$. It can be shown that
\bn\begin{array}{ll}
T^t&=v(uu_x+3tu_t+2u) ,\\
  T^x&=(uu_x+3tu_t+2u)(-uv-v_{xx})+v_x(3u_x+3tu_{xt}+xu_{xx})-v(4 u_{xx}+3tu_{xxt}+xu_{xxx}  )            \\
\end{array}\label{h}\en
so that
\bn\begin{array}{ll}
D_tT^t+D_xT^x&=(uu_x+3tu_t+2u  )( v_t-v_{xxx}-uv_x  )\\
&+  v5(u_t-uu_x-u_{xxx}  )+v(3tD_t+xD_x  )(u_t-uu_x-u_{xxx}  )   \\
&=W(v_t-v_{xxx}-uv_x )+v(3+  \R)  (u_t-uu_x-u_{xxx}  ),
\end{array}\label{h1}\en
where $\R=3tD_t+xD_x+2$, $W=uu_x+3tu_t+2u$ and   we note that $X(u_t-uu_x-u_{xxx} )=(3+  \R)  (u_t-uu_x-u_{xxx}  )$.

{\bf 4. Discussion}

It is clear that in each case, the conserved flows $(T^t,T^x)$ are nontrivial since $D_tT^t+D_xT^x$ do not vanish identically but, rather, on the solutions of the differential equation. The dependence of this method on solutions of the adjoint equation is equivalent to the multiplier approach since multipliers are solutions of the adjoint equation. Thus, as mentioned before, the two approaches in \cite{bgi} and \cite{i} complement each other and the latter has a formal procedure to construct the conserved flows
using symmetries of the differential equation. Moreover, we showed that the total divergence, quite explicitly, displays a relationship between symmetries (point or generalised) and conservation laws in a general setting - compare this to the results in \cite{ukm}. Also, the main results of this paper  mimics, to some extent, the results established on the relationship between symmetries and  multipliers of a differential equation as discussed in \cite{AK}.

\end{document}